\documentclass[10pt]{article}
\usepackage{amsfonts,amsmath,enumerate,amsthm,amssymb,color}
\numberwithin{equation}{section}
\theoremstyle{plain}
\newtheorem{thm}{Theorem}[section]

\newtheorem{cor}[thm]{Corollary}
\newtheorem{prop}[thm]{Proposition}

\theoremstyle{definition}

\usepackage{color}
\newcommand{\al}{\alpha}
\newcommand{\bt}{\beta}
\newcommand{\gm}{\gamma}
\newcommand{\Gm}{\Gamma}
\newcommand{\dl}{\delta}
\newcommand{\Dl}{\Delta}
\newcommand{\ep}{\varepsilon}

\newcommand{\Ld}{\Lambda}
\newcommand{\sg}{\sigma}

\newcommand{\ve}{\varepsilon}

\newcommand{\q}{\quad}

\newcommand{\R}{\mathbb{R}}
\newcommand{\N}{\mathbb{N}}

\newcommand{\law}{\mathcal L}

\newcommand{\bp}{\boxplus}

\newcommand{\intR}{\int_{\R}}
\newcommand{\sumn}{\sum_{j=1}^n}

\begin{document}
\title{Rates of convergence in the free central limit theorem}
\author{Makoto Maejima\\
Keio University 
\and Noriyoshi Sakuma\footnote{NS is supported by JSPS Kakenhi 19H01791, 19K03515, JPJSBP120209921, JPJSBP120203202.}\\Graduate School of Natural Sciences, \\Nagoya City University
}
\maketitle
\begin{abstract}
We study the free central limit theorem for not necessarily identically distributed free random variables where the limiting distribution is the semicircle distribution.
Starting from an estimate for the Kolmogorov distance between the measure of suitably normalized sums of free random variables and the semicircle distribution without any moment condition, we show the free Lindeberg central limit theorem and improve the known results on rates of convergence under the conditions of the existence of the third moments.
\end{abstract}


\vskip 10mm
\section{Introduction}

In classical probability theory, the central limit theorem is undoubtedly one of the most important problems.
Therefore, the situation is similar in non-commutative probability theory.
In this paper, we treat the central limit theorem in free probability theory, called the free central limit theorem.
Although the limiting distribution in the classical central limit theorem is Gaussian distribution, it is known that the limiting distribution in the free central limit theorem is the semicircle distribution. 
Due to this fact, the semicircle distribution is sometimes called free Gaussian distribution.
As for the origin of the free central limit theorem, see, among others, \cite{V85} and Theorem 3.5.2 in \cite{VDN92}.

With respect to the rate of convergence, Kargin \cite{Kargin07} and Chistyakov and G{\"o}tze \cite{CG08} proved Berry-Esseen type theorems for the free central limit theorem in the scalar-valued free probability setting via complex analysis method.
Mai and Speicher \cite{MS13} also proved it in the multivariate and operator-valued free probability setting. 
Some recent achievements on the Berry-Esseen type theorems in  free probability theory are
\cite{Austern} and \cite{BM}.
Austern \cite{Austern} gave a dynamical version of it, and Banna and Mai \cite{BM} improved its quantitative estimate and proved the Berry-Esseen type theorem in terms of the L{\'e}vy distance in the operator-valued free probability setting. All of their results used the operator-valued Cauchy transform technique. 
Among the above, only Chistyakov and G{\"o}tze \cite{CG08} treated not necessarily identically distributed free random variables,
which we state in Proposition~\ref{CG} below and relies on the existence of the third moments.
The purpose of this paper is to improve their result under weaker moment conditions.
Our approach does not use the complex analysis method but the method of truncation, which is popular in classical probability theory.


We now formulate the problem.
Let $\mathcal P(\R)$ be the set of all probability measures on $\R$.
For $\mu\in \mathcal P (\R)$ and $k\in\N$, define $m_k(\mu) = \intR x^k\mu(dx)$ and $\gm_k(\mu)= \intR |x|^k\mu(dx)$.
$\mu_w\in\mathcal P(\R)$ denotes the standard semicircle distribution with mean $0$ and variance $1$, having the probability distribution
\begin{align}\label{semic}
\mu_w (dx) = \frac{1}{2\pi}\sqrt{4-x^2}{\bf 1}_{\{|x|\le 2\}}(x)dx.
\end{align}
Define by $\Dl (\mu, \nu)$ the Kolmogorov distance between $\mu\in\mathcal P(\R)$
and $\nu\in\mathcal P(\R)$:
$$
\Dl(\mu, \nu) = \sup_{x\in\R}|\mu((-\infty, x]) - \nu((-\infty, x])|.
$$

Let $n\in\N$.
For $j=1,2, ..., n$, let $\mu_j\in\mathcal P(\R)$ and suppose that $m_1(\mu_j)=0$,
$\sg_j^2 = m_2(\mu_j)<\infty$ and $B_n^2= \sumn\sg_j^2>0$.
Write
$$
\mu_{n,j}((-\infty, x]) = \mu_j ((-\infty, B_nx]), \q x\in\R,\q j=1,2, ..., n.
$$


Denote by $\law(X)$ the probability measure of (classical or free) random variable $X$.
In classical probability theory, for $\mu, \nu \in \mathcal P(\R)$, 
$\mu*\nu$ denotes the convolution of $\mu$ and $\nu$, namely, 
$\mu*\nu = \law(X+Y)$,
where $X$ and $Y$ are independent random variables with $\mu= \law(X)$
and $\nu=\law(Y)$.
The free analogue of $\mu*\nu$ is written as $\mu\bp\nu$, called
the free additive convolution of $\mu$ and $\nu$,
which is $\mu\bp\nu = \law(X+Y)$, where
$X$ and $Y$ are free independent random variables with $\mu= \law(X)$
and $\nu=\law(Y)$.
(See, for instance, \cite{BV93}.)

In what follows, we use the following notation.
$$
\bp_{j=1}^n\mu_j = \mu_1\bp \cdots \bp \mu_n \q {\rm for}\q \mu_1,\ldots ,\mu_n\in \mathcal P (\R).
$$
Write $\mu^{(n)}= \bp_{j=1}^n\mu_{n,j}.$
The following rate of convergence in the free central limit theorem was shown in \cite{CG08}.
Throughout this paper, $C$ denotes a positive constant which may differ from one equality or inequality to another.

\begin{prop} [Theorem 2.6 in \cite{CG08}]\label{CG}
Let $n\in\N$ and suppose $\gm_3(\mu_j)<\infty$ for all $j=1,2,... ,n$.
Then there exists $C>0$ such that 
\begin{eqnarray}\label{1}
\Dl(\mu^{(n)}, \mu_w) \le \frac{C}{B_n^{3/2}}\left(\sumn\gm_3(\mu_j)\right)^{1/2}.
\end{eqnarray}
\end{prop}

Our aim of this paper is to relax the moment condition $\gamma_3(\mu_j)<\infty$ to
weaker moment conditions. 
We first, in Section 2, show an upper bound for the distance between $\bp_{j=1}^n\mu_j$
with a suitable normalization and $\mu_w$ 
without any moment condition, which is a free version of Theorem 7 in Chapter V of \cite{P75}.
In Section 3, we show an upper bound for $\Dl(\mu^{(n)}, \mu_w)$ under only the condition of the existence of the second moments.
In Section 4, as applications of the result shown in Section 3,  we show a free version of Lindeberg central limit theorem and 
the aim of this paper, namely, rates of convergence
under weaker moment conditions than the existence of the third moments.


\section{An estimate without any moment condition}


Let $n\in\N$ and let 
$\mu_j, j=1,2,...,n$, be probability measures on $\R$.
For $j=1,2, ... , n$, let $t_j < 0 < \tau_j$, $I_j=[t_j, \tau_j]$,
$$
\nu_j((-\infty, x]) =\left\{
\begin{array}{lll}
0 & (x<t_j)\\
\mu_j([t_{j}, x]) +\mu_j(I_j^{c})\delta_{0}((-\infty, x]) & (|x|\in I_j)\\
1 &  (x>\tau_j),
\end{array}
\right.
$$

$$
\al_j=\int_{\R}x\nu_j(dx) \q{\rm  and}\q
\bt_j^2 = \int_{\R}x^2\nu_j(dx) -\al_j^2.
$$
Assume that  $\bt^2\ne 0$ for some $j=1,2, ... , n$.
Let
$M_n =\sumn \al_j, N_n^2= \sumn \bt^2_j$
and
$$
\nu_{n,j}((-\infty, x]) = \nu_j((-\infty, N_nx + M_n ]),
\text{ namely, }
\nu_j((-\infty, x]) = \nu_{n,j}\left(\left(-\infty, \frac{x-M_n}{N_n}\right]\right).
$$
Let $\mu_w$ be the standard semicircle distribution defined in \eqref{semic} and let
$$
\nu^{(n)} = \bp_{j=1}^n \nu_{n,j},\q
\Dl_n = \Dl(\nu^{(n)}, \mu_w),\q
\Gm_n = \sumn \mu_j(I_j^c).
$$

\begin{thm}\label{nomoment}
For any $a>0$ and $b\in\R$,
\begin{align*}
\sup_{x\in\R}& \left| (\bp_{j=1}^n\mu_j)((-\infty, a(x+b)]) - \mu_w((-\infty, x])\right|
\le \Dl_n + \Gm_n + C\frac{|ab -M_n|}{N_n} + C\left|\frac{a}{N_n}-1\right| .
\end{align*}
\end{thm}

{\it Proof.}
It follows from the defintion of $\nu_j$ that
$$
|\mu_j((-\infty, x]) - \nu_j((-\infty,x])| \le \mu_j(I_j^c)\q {\rm for \,\,every}\,\, x\in\R,
$$
implying that
$$
\Dl(\mu_j, \nu_j) \le \mu_j (I_j^c).
$$
Proposition 4.13 in \cite{BV93} says that
if $\mu, \mu', \nu, \nu'\in \mathcal P(\R)$, then
$$
\Dl(\mu\bp\nu, \mu'\bp \nu') \le \Dl(\mu, \mu') + \Dl(\nu, \nu').
$$
Thus by the repeated use of this result, we have
$$
\Dl\left( \bp_{j=1}^n \mu_j , \bp_{j=1}^n \nu_j\right)
\le \sumn \Dl(\mu_j, \nu_j)\le \sumn \mu_j(I_j^c) =\Gm_n.
$$
We now have
\begin{align*}
\Big| (\bp_{j=1}^n \mu_j) & ((-\infty, a(x+b)]) - \mu_w((-\infty, x])\Big|\\
&\le 
\left|(\bp_{j=1}^n \mu_j)((-\infty, a(x+b)]) - (\bp_{j=1}^n \nu_j)((-\infty, a(x+b)])\right| \\
&\hskip 10mm+ \left| (\bp_{j=1}^n \nu_j)((-\infty, a(x+b)])-\mu_w((-\infty, px+q])\right| \\
&\hskip 10mm+ |\mu_w((-\infty, px+q]) - \mu_w((-\infty, px]) 
+ |\mu_w((-\infty, px]) - \mu_w((-\infty, x])|
\end{align*}
for every $p,q\in\R$.
Let
$$
p=\frac{a}{N_n},\q
q=\frac{ab-M_n}{N_n}.
$$
Then
\begin{align*}
\Big| (\bp_{j=1}^n
 \nu_j)((-\infty, a(x+b)]) -\mu_w((-\infty, px+q])\Big|
= \left| \nu^{(n)}((-\infty, px +q])
-\mu_w((-\infty, px + q])\right| \le \Dl_n.
\end{align*}
Therefore, we have
$$
\sup_{x\in\R} \left|(\bp_{j=1}^n \mu_j)((-\infty, a(x+b)]) - \mu_w((-\infty, x])\right|
\le \Gm_n +\Dl_n +  T_1 + T_2,
$$
where
$$
T_1= \sup_{x\in\R} |\mu_w((-\infty, x+q]) -\mu_w((-\infty, x])|
\q\text{ and }\q
T_2= \sup_{x\in\R} |\mu_w((-\infty, px]) -\mu_w((-\infty, x])|.
$$
It is easy to see that, for $q\in\R$,
\begin{align}\label{w-ineqplus}
\sup_{x\in\R}|\mu_w((-\infty, x+q]) - \mu_w((-\infty, x])| \le \frac1{\pi}|q|.
\end{align}
and for $p>0$,
\begin{align}\label{w-ineqx}
\sup_{x\in\R}| \mu_w((-\infty, px]) - \mu_w((-\infty, x])|
\le \sup_{x\in\R} \frac1{2\pi}\left| \int_x^{px} \sqrt{4-u^2}{\bf 1}_{\{|u|\le 2\}}(u)du\right|
\le \frac2{\pi}|p-1|
\end{align}
The inequality \eqref{w-ineqplus} implies that
$$
T_1 \le \frac1{2\pi}|q| = C\left|\frac{ab-M_n}{N_n}\right|
$$
and
\eqref{w-ineqx} implies that
$$
T_2\le \frac{4}{2\pi}|p-1| = C\left|\frac{a}{N_n}-1\right|.
$$
The proof of Theorem~\ref{nomoment} is now completed.


\section{An estimate under the condition of the existence of second moments}

By the defintion of $\nu_j$, all measures $\nu_j$ have bounded supports, implying that they have all finite moments.
Therefore we can apply Proposition~\ref{CG} to the estimate of $\Dl_n$ in Theorem~\ref{nomoment}.
Then we have
\begin{align}\label{Delta}
\Dl_n \le \frac{C}{N_n^{3/2}}\left(\sumn \int_{\R} |x-\al_j|^3\nu_j(dx)\right)^{1/2}
\le \frac{{2\sqrt{2}} C}{N_n^{3/2}}\left(\sumn \int_{\R} |x|^3\nu_j(dx)\right)^{1/2},
\end{align}
where the last inequality holds by the fact that $|a+b|^{p} \le 2^{p-1}(|a|^{p} + |b|^{p})$ for $p\ge 1$.
\eqref{Delta} will be used in the proof below.

\begin{thm}\label{no-moment}
For $j=1,2, ... ,n$, let $\mu_j\in\mathcal P(\R)$ and suppose $m_1(\mu_j)=\int_{\R}x\mu_j(dx)=0$, $\sg_j^2 =\int_{\R}x^2\mu_j(dx)<\infty$ and $B_n^2=\sumn \sg_j^2>0$.
For $\ep \in (0,1]$, let
$$
\Ld_n(\ve) = \frac1{B_n^2}\sumn \int_{|x|> \ve B_n}x^2\mu_j(dx)
\q\text{ and }\q
\ell_n(\ve) = \frac1{B_n^3}\sumn\int_{|x|\le \ve B_n}|x|^3\mu_j(dx).
$$
Then
\begin{align}\label{epbound}
\Dl(\mu^{(n)}, \mu_w) \le C (\Ld_n(\ve) + \ell_n(\ve))^{1/2} \q {\rm for\,\, every}\,\,\ve\in(0,1].
\end{align}
\end{thm}

\vskip 3mm
{\it Proof}.
We first prove \eqref{epbound} for $\ep=1$.
In Theorem~\ref{nomoment}, we let $a=B_n, b=0$,
$-t_j =\tau_j =B_n$ for $j=1,2, ... ,n.$
Note that $N_n^2\le B_n^2$.
If $N_n^2\le \frac14B_n^2$, then
\begin{align}
\frac34B_n^2 & \le B_n^2 - N_n^2
 = \sumn \left\{\int_{|x|> B_n}x^2\mu_j(dx) + \left(\int_{|x|\le B_n}x\mu_j(dx)\right)^2\right\}
 \le 2 \sumn \int_{|x| > B_n}x^2 \mu_j(dx),
\end{align}
and thus $\Ld_n(1) \ge \frac38$.
Since $\Dl(\mu^{(n)}, \mu_w) \le 1$, \eqref{epbound} holds with $C=\frac83$.

We next consider the case when $N_n^2> \frac14B_n^2.$
It follows from \eqref{Delta} that
\begin{align*}
\Dl_n & \le \frac{2\sqrt{2}C}{N_n^{3/2}}\left(\sumn \int_{|x|\le B_n}|x|^3\mu_j(dx)\right)^{1/2}\le \frac{8C}{B_n^{3/2}}\left(\sumn \int_{|x|\le B_n}|x|^3\mu_j(dx)\right)^{1/2} = C\ell_n(1)^{1/2}.
\end{align*}
Furthermore, 
$$
\Gm_n = \sumn \int_{|x|>B_n} \mu_j(dx)
\le \frac1{B_n^2}\sumn \int_{|x| > B_n} x^2\mu_j(dx)
= \Ld_n(1),
$$
$$
\frac{|M_n|}{N_n}  \le \frac2{B_n}\sumn \int_{|x| >  B_n} |x|\mu_j(dx)\le2 \Ld_n(1)
$$
and
\begin{align*}
\left|\frac{B_n}{N_n} -1\right| & \le \left|\frac{B_n}{N_n} -1\right|\cdot
\left|\frac{B_n}{N_n} +1\right|\le \frac4{B_n^2}|B_n^2-N_n^2|= \frac4{B_n^2}\sumn \int_{|x| > B_n}x^2\mu_j(dx) 
 = 4\Ld_n(1).
\end{align*}
Thus, by combining above estimates and by Theorem~\ref{nomoment}, we have
$$
\Dl(\mu^{(n)}, \mu_w) \le C(\Ld_n(1) + \ell_n(1)^{1/2}),
$$
and \eqref{epbound} holds for $\ep = 1$.

Now let $\ep\in (0,1]$ arbitrarily.
Then $\Ld_n(1)\le \Ld_n(\ep) \le \Ld_n(\ep)^{1/2}$ (since $\Ld_n(\ep)\le 1$),  and
$$
\ell_n(1) = \ell_n(\ve) + \frac1{B_n^3}\sumn \int_{\ep B_n < |x| \le B_n}|x|^3\mu_j(dx)\le
\ell_n(\ve) + 
\frac1{B_n^2}\int_{|x|> \ep B_n}x^2 \mu_j(dx) \le
\ell_n(\ep) + \Ld _n(\ep).
$$
Consequently, we have
$$
\Ld_n (1) + \ell_n(1)^{1/2} \le \Ld_n(\ep)^{1/2} +  (\Ld_n(\ep) + \ell_n (\ep))^{1/2}
\le 2 (\Ld_n(\ep) + \ell_n (\ep))^{1/2}
$$
for every $\ep\in (0,1]$.
This completes the proof.


\section{The free Lindeberg central limit theorem and the Berry-Esseen type estimates.}

In this section, we first show, as an application of Theorem~\ref{no-moment}, the following theorem, which is a free version of a well-known result in classical probability theory. The condition \eqref{Lind} below is called the 
Lindeberg condition.

\begin{thm}[Free Lindeberg central limit theorem]
If for every fixed $\ep>0$,
\begin{align}\label{Lind}
\Ld_n(\ve)= \frac1{B_n^2}\sumn \int_{|x| >  \ve B_n}x^2\mu_j(dx)\to 0 \q {\rm as}\,\, n\to \infty,
\end{align}
then the free central limit theorem
holds, namely,
$$
\Dl (\mu^{(n)}, \mu_w) \to 0 \q {\rm as}\,\, n\to \infty.
$$
\end{thm}

{\it Proof.}
From the definition of $\ell_n(\ep)$, we see that
$$
\ell_n(\ep)\le \frac{\ep}{B_n^2}\sumn \int_{|x|\le \ep B_n}x^2 \mu_j (dx) \le \ep.
$$
Hence, by \eqref{epbound}, we have
$
\Dl(\mu^{(n)}, \mu_w) \le C( \Ld_n(\ep)+\ep)^{1/2}
$
for every $\ep\in (0,1]$.
This concludes the statement.

\vskip 3mm

We finally show rates of convergence in the free central limit theorem
under weaker moment conditions than the existence of the third moments.
, by using Theorem~\ref{no-moment}.

Let $G$ be the set of real-valued functions, defined on $\R$, that satisfy the conditions:
(a) $g(x)$ is nonnegative, even and nondecreasing on $(0,\infty)$ and 
(b) $\frac{x}{g(x)}$ is nondecreasing on $(0,\infty)$.
A typical example is $g(x) = |x|^{\dl}$,  $0<\dl\le 1$.
Write $\gm_{2+g}(\mu_j) = \intR x^2g(x)\mu_j(dx)$.

\begin{thm}\label{main}
Let $g\in G$ and suppose $\gm_{2+g}(\mu_j)<\infty$ for $j=1,2,... ,n$.
Then there exists $C>0$ such that
\begin{eqnarray}\label{2}
\Dl(\mu^{(n)}, \mu_w) \le \frac{C}{B_ng(B_n)^{1/2}}\left(\sumn \gm_{2+g}(\mu_j)\right)^{1/2}.
\end{eqnarray}
\end{thm}

{\it Proof.}
We have
\begin{align*}
\Ld_n(1) &= \frac1{B_n^2}\sumn \int_{|x|> B_n}x^2\mu(dx) = \frac1{B_n^2g(B_n)}\sumn 
\int_{|x|> B_n} x^2g(x) \frac{g(B_n)}{g(x)}\mu_j(dx).
\end{align*}
Since $g$ is nondecreasing on $(0,\infty)$,
$$
\frac{g(B_n)}{g(x)}\le 1 \q{\rm on}\,\, \{x: |x|>B_n\}.
$$
Hence 
$$
\Ld_n(1) \le \frac{1}{B_n^2g(B_n)}\sumn \gm_{2+g}(\mu_j).
$$
We also have
\begin{align*}
\ell_n(1)& = \frac1{B_n^3}\sumn\int_{|x|\le B_n}|x|^3\mu_j(dx) = \frac1{B_n^2g(B_n)}\sumn 
\int_{|x|\le B_n}x^2g(x) \frac{\left(\frac{|x|}{g(x)}\right)}{\left(\frac{B_n}{g(B_n)}\right)}\mu_j(dx).
\end{align*}
Since $\frac{x}{g(x)}$ is nondecreasing on $(0,\infty)$,
$$
\frac{\left(\frac{|x|}{g(x)}\right)}{\left(\frac{B_n}{g(B_n)}\right)}\le 1\q{\rm on}\,\, \{x: |x|\le B_n\}.
$$
Hence
$$
\ell_n(1) \le \frac{1}{B_n^2g(B_n)}\sumn \gm_{2+g}(\mu_j).
$$
Then Theorem~\ref{no-moment} implies that 
$$
\Dl(\mu^{(n)}, \mu_w) \le C(\Ld_n (1) + \ell_n(1) )^{1/2}\le \frac{C}{B_ng(B_n)^{1/2}}
\left(\sumn \gm_{2+g}(\mu_j)\right)^{1/2},
$$
completing the proof.

\vskip 3mm
If we choose $g(x)=|x|^{\dl}, 0<\dl\le 1$, we have the following.

\begin{cor}\label{cor}
Let $0<\dl \le 1$ and suppose $\gm_{2+\dl}(\mu_j)=\int_{\R}|x|^{2+\dl}\mu_j(dx)<\infty$ for $j=1,2,,... ,n$.
Then there exists $C>0$ such that
\begin{eqnarray}\label{3}
\Dl(\mu^{(n)}, \mu_w) \le \frac{C}{B_n^{1+\dl/2}}\left(\sumn \gm_{2+\dl}(\mu_j)\right)^{1/2}.
\end{eqnarray}
\end{cor}

When $\dl=1$, \eqref{3} is \eqref{1} in Proposition~\ref{CG}.


\section*{Acknowledgements}
The authors would like to thank Leonie Neufeld for
their helpful comments.
\vskip 5mm

\end{document}